\documentclass[11pt]{article}
\usepackage[margin=4.46cm]{geometry}
\usepackage{graphicx}
\usepackage{epsfig}
\usepackage{amsmath}
\usepackage{amssymb}
\usepackage{color}
\usepackage{caption}
\usepackage{tikz}
\usetikzlibrary{decorations.markings}

\newcommand{\commentout}[1]{}

\def \Rset {{\mathbb R}}

\newcommand{\nit}{\noindent}

\newcommand{\be}{\begin{equation}}
\newcommand{\ee}{\end{equation}}
\newcommand{\ba}{\begin{eqnarray}}
\newcommand{\ea}{\end{eqnarray}}
\newcommand{\bi}{\begin{itemize}}
\newcommand{\ei}{\end{itemize}}
\newcommand{\br}{\begin{eqnarray}}
\newcommand{\er}{\end{eqnarray}}

\newcommand{\qed}{\mbox{$\square$}\newline}

\newtheorem{theo}{Theorem}[section]

\newtheorem{lem}{Lemma}[section]
\newtheorem{cor}{Corollary}[section]

\begin{document}
\title{Periodic orbits of the ABC flow with $A=B=C=1$ }
\author{ Jack Xin\footnote{Department of Mathematics, UC Irvine, jxin@math.uci.edu, partially supported 
by NSF grant DMS-1211179.},\quad  Yifeng Yu\footnote{Department of Mathematics,  UC Irvine, yyu1@math.uci.edu,  
partially supported by NSF CAREER grant DMS-1151919.},\quad Andrej Zlato\v{s}\footnote{Department of Mathematics, 
UW--Madison, andrej@math.wisc.edu,  partially supported by NSF CAREER grant DMS-1056327.}}

\date{}
\maketitle

\begin{abstract}
\nit In this paper,  we prove that the ODE system
%\begin{cases}
\begin{align*}
\dot x &=\sin z+\cos y\\
\dot y &= \sin x+\cos z\\
\dot z &=\sin y + \cos x,
\end{align*}
%\end{cases}
whose right-hand side is the Arnold-Beltrami-Childress (ABC) flow with parameters $A=B=C=1$,
has periodic orbits on $(2\pi\mathbb T)^3$
%(when viewed on $(2\pi\mathbb T)^3$) 
with rotation vectors parallel to $(1,0,0)$, $(0,1,0)$, and $(0,0,1)$.   
An application of this result is  that the well-known G-equation model for turbulent combustion with this ABC flow on $\mathbb R^3$ has a linear (i.e., maximal possible) flame speed enhancement  rate as the amplitude of the flow grows.
\end{abstract}

\section{Introduction}  
 Arnold-Beltrami-Childress (ABC) flows are important examples of 
three dimensional periodic incompressible flows \cite{Ar65,DFGHMS1986,Fr93}. The general form is 
\[
u_{A,B,C}(x,y,z)=(A\, \sin z+  C \, \cos y, B\, \sin x + A\, \cos z, C \,  \sin y +  B\, \cos x )
\]
%\be\label{abc}
%\begin{split}
%\dot x &= A\, \sin z+ \,  C \, \cos y \\
%\dot y &=  B\, \sin x + \,  A\, \cos z \\
%\dot z &= B\, \cos x + \,  C \,  \sin y,
%\end{split}
%\ee
with $A,B,C\in\Rset$ some given parameters.   These flows are steady solutions of the Euler equation, 
and the corresponding system of ODE
\[
\dot X(t)=u_{A,B,C}(X(t))
\]
becomes integrable when one of the parameters is zero.
We refer the reader to \cite{DFGHMS1986} and to Section~2.6 in \cite{C1995} for more background information.

The most interesting case is when $A=B=C=1$ (the 1-1-1 ABC flow), 
introduced by Childress \cite{C1967,C1970} in connection with dynamo studies. 
The 1-1-1 ABC flow system is known to be chaotic in the sense that a web of chaos 
occupies a (fairly complex) region of 
the phase space, as seen from the Poincar\'e section plots of Figure 10 in 
\cite{DFGHMS1986}. Though chaotic orbits for other ABC parameters have been identified 
analytically (see, e.g.,  \cite{Ar65,zhkwlihu93}), it is not known 
how to construct them rigorously for the 1-1-1 ABC flow. 
%\medskip
 
In this paper, we prove that the 1-1-1 ABC flow system has unbounded  orbits which are periodic up to shifts by $2\pi\mathbb Z$-multiples of one of the standard basis vectors. 
Our interest in this question is motivated by recent progress as well as open questions in the study of flame speed
enhancement by turbulent fluid motion in combustion models 
%of front propagation of reaction-advection-diffusion partial differential equations 
(see, e.g., \cite{XY2014,Z_10}), especially in three spatial dimensions.  Indeed, linear (and hence maximal possible) growth of flame front speeds as the amplitude of the flow grows is related to
existence of unbounded (roughly) ``periodic'' orbits of the turbulent flow 
%is related to 
%linear growth of flame front speeds as the amplitude of the flow grows 
(see below), rather than to existence of  chaotic orbits (which one might expect). 
The study of the former orbits  thus becomes natural, and is a distinct 
departure from the long history of the study of chaos in ABC flows. 

%??? The traditional view is that chaotic streamlines (Lagrangian turbulence) may considerably enhance 
%transport. 
Note that since the 1-1-1 ABC flow is not a small perturbation of an integrable case, 
% and degenerate so that 
classical tools in dynamical systems (such as Melnikov analysis or Smale horseshoe) 
are difficult to apply. 
%Surprisingly, there is order within chaos that permits a novel mathematical construction of unbounded (extended) trajectories. 
Nevertheless, using the analysis of the orbits of the flow 
inside a triangular prism region (see Figure 1) and certain symmetries of the flow,
% play  a critical role in our proof. Below is our main theorem.
we are able to establish the following main result of this paper.

\begin{theo}\label{periodic}   There exists  $t_0>0$ and  a solution  $X(t)=(x(t),y(t), z(t))$ to the 1-1-1 ABC flow system such that for each $t\in\Rset$ we have
$$
X(t+t_0)=X(t)+(2\pi, 0, 0).
$$  
Then  $Y(t)=(z(t), x(t), y(t))$ and $Z(t)=(y(t),   z(t),  x(t))$ are clearly also solutions and satisfy 
$$
Y(t+t_0)=Y(t)+(0, 2\pi, 0)   \qquad \mathrm{and}   \qquad Z(t+t_0)=Z(t)+(0, 0, 2\pi)  .
$$
\end{theo}

\medskip

%???ALSO $A,B,C\sim 1$.
%
%\medskip

{\it Remarks.}
%A few remarks are in order.
%\nit $\bullet$ {\bf Other cases.}
1. Then $\tilde X(t)=X(-t)-(\pi,\pi,\pi)$, $\tilde Y(t)=Y(-t)-(\pi,\pi,\pi)$, and $\tilde Z(t)=Z(-t)-(\pi,\pi,\pi)$ are  also solutions and they satisfy 
$\tilde X(t+t_0)=\tilde X(t)-(2\pi, 0, 0)$, $\tilde Y(t+t_0)=\tilde Y(t)-(0, 2\pi, 0)$, and $\tilde Z(t+t_0)=\tilde Z(t)-(0, 0, 2\pi)$.\smallskip
  
2. Our method can be adjusted to obtain unbounded ``periodic'' orbits along the $x$, $y$ and $z$ directions for some other values of $A,B,C$. For instance,  if $A$, $B$ and $C$ are all close to 1, then this follows from the proof of our result via a simple perturbation argument.   Another example is the case $0<A\ll 1$ and $B=C=1$,  which is a perturbation of the integrable case $A=0$ and $B=C=1$.  We note that the corresponding system now possesses a large KAM region as well as a chaotic thin layer near the separatrix walls of the integrable flow (i.e., near  $\{\sin y+\cos x=0\}$), and we refer to \cite{MXYZ2015} for  both theoretical  and numerical analysis of this case.  

\medskip

\nit {\bf  Applications to combustion models.}   
%Our motivation to finding  unbounded periodic orbits comes from the study of flame propagation. 
Finding the turbulent flame speed (or effective burning velocity) is one of the most important unsolved problems in turbulent combustion.  Roughly speaking, turbulent flame speed is the flame propagation speed in the presence of (and enhanced by)  a strong flow of the ambient fluid medium.  Two typical examples are the spread of wildfires fanned by winds and combustion  of rotating air-gasoline mixtures inside internal combustion engines.   

For simplicity,  let us assume that the flow velocity profile $V:\Rset^n\to \Rset^n$ is smooth, periodic, and incompressible (i.e., $\nabla\cdot V=0$).  A well-known approach to the study of flame propagation and turbulent flame speed is the G-equation model (see, e.g., \cite{Pet00,W85}),  the level set Hamilton-Jacobi equation
$$
G_t+KV(x)\cdot \nabla G+s_l|\nabla G|=0,
$$
with $K\ge 0$ the amplitude of the turbulent flow and $s_l>0$ the laminar flame speed.  A basic question is to understand how the turbulent flame speed depends on the flow amplitude as $K\to\infty$ (i.e., for strong flows).  Let $s_T(p,K)$ be the turbulent flame speed given by the G-equation model along a fixed unit direction $p\in\Rset^n$
% when $s_l$ is a constant.  
(see  \cite{XY2014} for the precise definition and  further references).  Then  it  is proved in \cite{XY2014} that 
$$
\lim_{K\to \infty}{s_T(p,K)\over K}=\max_{\{\xi|\ \dot \xi=V(\xi)\}}\limsup_{t\to \infty}{p\cdot \xi(t)\over t }.
$$
In particular,  $s_T(p,K)$ grows linearly  (which is the maximal possible growth rate) as $K\to\infty$ precisely when 
there exists an orbit of $\dot \xi=V(\xi)$ which travels roughly linearly in the direction $p$.

This yields the following corollary of  Theorem \ref{periodic} (and of Remark 1):

\begin{cor}  If $V$ is the 1-1-1 ABC flow, then for any  $p\in\Rset^3$ we have 
$$
\lim_{K\to \infty}{s_T(p,K)\over K}>0.
$$
\end{cor}

Another well-known model used in the study of turbulent flame speeds involves traveling front solutions of the reaction-advection-diffusion equation
$$
T_t+KV(x)\cdot \nabla T=d\Delta T+f(T).
$$
Here $T$ represents the  temperature of the reactant, $d>0$ is the molecular diffusivity, and $f$ 
is a nonlinear reaction function (see, e.g., \cite{Berrev, Xin_09}).   Consider the case of a KPP reaction $f$ (e.g., $f(T)=T(1-T)$) 
and let $c^{*}(p,K)$ be the turbulent flame speed in the direction $p$ given by this model (i.e., the minimal speed of a traveling front in direction $p$; see \cite{Berrev, Xin_09} for details).  It is established in  \cite{Z_10} that
$$
\lim_{K\to \infty}{c^{*}(p,K)\over K}=\sup_{w\in\Gamma}\int_{\Bbb T^n}(V\cdot p) w^2dx,
$$
%where $\sigma=w^2dx$ for $w\in  \Gamma$ defined  as ($H^1$ invariant  probability measures)
where
$$
\Gamma= \left \{w\in  H^1(\Bbb T^n) \,\Big|\,   V\cdot \nabla w=0 \,\&\, ||w||_{L^2(\Bbb T^n)}=1 \,\&\, \|\nabla w||_{L^2(\Bbb T^n)}^{2}\leq f'(0) \right\}.
$$
Hence,  in contrast to the G-equation model, one now needs a positive measure of orbits of $\dot \xi=V(\xi)$ which travel roughly linearly in the direction $p$ to obtain linear-in-$K$ turbulent flame speed enhancement.  

When $n=2$,  stability of periodic orbits was used in \cite{XY2014} to establish   
that  $ \lim_{K\to \infty}{c^{*}(p,K)\over K}=0$ if and 
only $\lim_{K\to \infty}{s_T(p,K)\over K}=0$.  However, this is not true in general in three dimensions.  An example from  \cite{XY2014}  is the so-called  Robert cell flow, for which $\lim_{K\to \infty}{s_T(p,K)\over K}>0$ and $\lim_{K\to \infty}{c^{*}(p,K)\over K}=0$ when $p=(0,0,1)$.   

The analysis becomes much more difficult for the  much more interesting 1-1-1 ABC flow  due to the presence of chaotic structures.   Nevertheless, numerical simulations suggest that there are ``vortex tubes'' composed of  orbits which travel roughly linearly in the $\pm x$, $\pm y$, and $\pm z$ directions \cite{DFGHMS1986}, suggesting that
$c^{*}(p,K)$ should also grow linearly in $K$ for each $p$.   
%However, we do not know of a rigorous proof of 
%this.  In particular,  we do not know whether the ``periodic'' orbits constructed in this paper are isolated or whether there are 
%other (quasi-)periodic orbits near them.
One could expect to find such vortex tubes, if they indeed exist, in the vicinity of the ``periodic'' orbits constructed in this paper, but  a rigorous proof of their existence is currently not known.

%%%%%%%%%%%%%%%%%%%%%%%%%%%%%%%%
\section{Proof of Theorem \ref{periodic}}  
Let $R$ be the open triangle in the $xy$-plane with vertices $(0, -{\pi \over 2})$, $(0,{3\pi\over 2})$, $(-\pi, { \pi \over 2})$, and let  $D=R\times(0,\frac\pi 2)$.  Our proof is based on showing that there exists a solution $X_{\bar a}$ to the 1-1-1 ABC flow system which starts from $(-{\pi\over 2}, 0,  \bar a)$ for some $\bar a\in  [0,  {\pi\over 2})$ and passes through the segment \hbox{$\{(0, y, {\pi\over 2})|\  y\in  [-{\pi\over 2},  {3\pi\over 2}]\}$} (see Figure 1).  We do this in Step 1 below, and then use the symmetries of the flow to construct the desired solution $X$ in Step 2.

\tikzset{->-/.style={decoration={
  markings,
  mark=at position .53 with {\arrow{>}}},postaction={decorate}}}

\begin{figure}[h] %\label{fig}
\begin{center}
\begin{tikzpicture}

%\fill (0,0) circle (0.08);

\draw[thick] (-1.5, -2)--(1.5,0)--(-4.5,0)--(-1.5,-2);
%\draw[dotted] (-4.5,-2.2)--(1.5,-2.2);
\draw[thick] (-4.5,-2.2)--(-1.5, -4.2)--(1.5,-2.2);

\draw[thick](-1.5, -2)--(-1.5, -4.2);
\draw[thick](1.5,0)--(1.5, -2.2);
\draw[thick](-4.5,0)--(-4.5,-2.2);
%\draw[thin, dotted](0,-3.2)--(0,-1);
\draw (-0.3,-2.4) node {$(-{\pi\over 2}, 0,  \bar a)$};
\fill (0,-2) circle (0.08);
\fill (-1,0) circle (0.08);
\draw[very thick, dotted,->-] (0,-2) to [out=90,in=195] (-1,0);
\draw(-1.6,-0.8) node {$X_{\bar a}(t)$};

\draw (-4.5,-2.8) node {$z=0$};
\draw (-4.5,0.4) node {$z={\pi\over 2}$};

\end{tikzpicture}
\captionof{figure}{The region $D$, rotated counter-clockwise by 90 degrees, and $X_{\bar a}$.}
\end{center}
\end{figure}
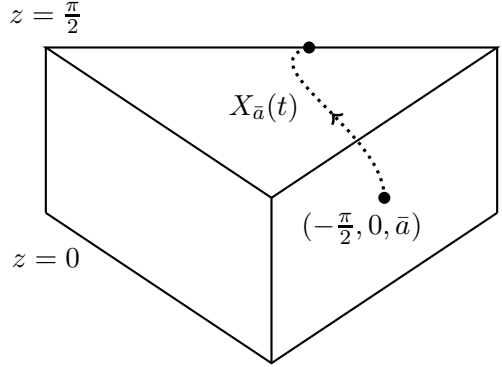

\medskip

\nit {\bf Step 1:}  
%We want to show that there exists a solution to the 1-1-1 ABC flow which starts from $(-{\pi\over 2}, 0,  \bar a)$ for some $\bar a\in  [0,  {\pi\over 2}]$ and passes through the line segment $ \{(0, y, {\pi\over 2})|\  y\in  [-{\pi\over 2},  {3\pi\over 2}]\}$.    
For any $a\in  [0,  {\pi\over 2})$,  let $X_a(t)=(x(t), y(t), z(t))$ satisfy $X_a(0)=(-{\pi\over 2}, 0, a)$ and
\begin{align*}
\dot x &=\sin z+\cos y\\
\dot y &=\sin x+\cos z\\
\dot z &=\sin y + \cos x.
\end{align*}
Obviously,   $(x(t),y(t))\in R$ for all small $t>0$ (and any $a\in  [0,  {\pi\over 2})$).  Since $\cos x +\sin y>0$ for $(x,y)\in R$, we have that  $X_a(t)\in D$ for all small $t>0$.   

The question now is whether and where $X_a$ will (first) exit $D$.  Clearly  $\dot z(t)>0$ when $X_a(t)\in D$, and we also have
$$
\dot X_a(t)\cdot (1,1,0)=\sin z(t)+\cos z(t)>0  % \quad \text{if $X_a(t)\in  \{x+y=-{\pi\over 2}\}\cap \partial D$}.
$$
when $X_a(t)\in  \{x+y=-{\pi\over 2}\}\cap \partial D$.
Hence  $X_a$ cannot exit $D$ through either the plane $\{z=0\}$ or  the plane $\{x+y=-{\pi\over 2}\}$. 

Let us now consider any $T>0$ such that $X_a(t) \in  \bar D=\bar R\times [0,  {\pi\over 2}]$ for all $t\in  [0, T]$.
The following result shows, in particular, that the first exit cannot happen through the plane $\{y-x={{3\pi}\over 2}\}$ either.  

\begin{lem}\label{l1}  
We have $y(t)< {\pi\over 4}$   for $t\in  [0,T]$.
%$$
%y(s)< {\pi\over 4}   \quad \text{for $s\in  [0,T]$}.
%$$
\end{lem}

{\it Proof.}  We argue by contradiction.  If not, then by $y(0)=0$, there are  $0\leq t_1<t_2\leq T $ such that $y(t)\in  [0,  {\pi\over 4}]$ for $t\in [t_1,t_2]$, $y(t_1)=0$, and $y(t_2)={\pi\over 4}$.  Let us choose the smallest such $t_1,t_2$.  Then $y(t)\in[-{\pi\over 2}, {\pi\over 4}]$ for  all $t\in[0,t_2]$ and therefore $\dot x(t)\ge 0$ for $t\in[0,t_2]$.
Therefore 
$$
w(t)=x(t)+{\pi\over 2}
$$
satisfies $w(0)=0$ and $w(t)\in [0,{\pi\over 2}]$ for  $t\in[0,t_2]$.

In the following consider only $t\in [t_1,t_2]$. We have   
$$
\dot z-\dot y=\sin y+\cos x-\sin x-\cos z\geq \sin y+1-\cos z\geq \sin y\geq 0,
$$
hence $z(t)\geq y(t)$ for $t\in [t_1,t_2]$ due to $z(t_1)\ge 0=y(t_1)$.   Then 
\begin{align*}
\dot w &=\sin z+\cos y\geq \sin y + \cos y\geq 1\\
\dot y &=\cos z-\cos w \leq 1-\cos w\leq \min \left\{{w^2\over 2},1\right\}
\end{align*}
since $w\in [0,{\pi\over 2}]$.  In particular, $\dot y\leq 1$ shows that  $t_2-t_1\geq {\pi\over 4}$.  Let us now consider three cases.

\medskip

Case 1:  $w(t_1)\geq {\pi\over 3}$.    Then $\dot w\geq 1$ shows  $w(t_2)\geq  {7\pi\over 12}>{\pi\over 2}$, a contradiction.

\medskip

Case 2:  $w(t_1)<{\pi\over 3}$ and $w(t_2)\leq {\pi\over 3}$.   We have $\dot y\leq  \dot w  {w^2\over 2}$, 
%Then
%$$
%\dot y\leq  {w^2\over 2}\leq \dot w  {w^2\over 2}.
%$$
%Taking integration on both sides over $[t_1,t_2]$,  we derive that  $y(t_2)\leq {w^3(t_2)\over 6}<0.2<{\pi\over 4}$.  A contradiction. 
which after integration over $[t_1,t_2]$ yields  $y(t_2)\leq {1\over 6} w^3(t_2)\le {{\pi^3}\over{162}}<{\pi\over 4}$, a contradiction.

\medskip

Case 3:   $w(t_1)<{\pi\over 3}$ and $w(t_2)> {\pi\over 3}$.   Then there exists $t^{*}\in  [t_1,t_2]$ such that  $w(t^{*})={\pi\over 3}$, and  the computation in Case 2 shows $y(t^{*})\le {{\pi^3}\over{162}}$.  Then 
\[
y(t_2)-y(t^{*})\ge {\pi\over 4}-{{\pi^3}\over{162}}>  {\pi\over 6}= {\pi\over 2}-{\pi\over 3} \ge w(t_2)-w(t^{*}),
\]
a contradiction with $\dot w\geq 1\geq \dot y$  on  $[t_1, t_2]$.  The proof is finished. \qed

\begin{lem}\label{l2}  
If $a=0$, then $x(t)+{\pi\over 2}>z(t)$ for all $t\in(0,T]$.
\end{lem}

{\it Proof.}  Again let $w(t)=x(t)+{\pi\over 2}$.  From $w(0)=y(0)=z(0)=0$  and 
\begin{align*}
\dot w & =\sin z+\cos y\\
\dot z & =\sin w+\sin y
\end{align*}
we have $\dot w(0)>\dot z(0)$,  hence $w(t)>z(t)$ for all small $t>0$.  Assume that there is $t^*\in (0,T]$ such that $w(t^*)=z(t^*)$ and $w>z$ on  $(0,t^*)$.  Lemma \ref{l1} shows that $\cos y>\sin y$ on $[0,T]$,  so
$$
0\geq \dot w(t^*)-\dot z(t^*)> \sin z(t^*)-\sin w(t^*)=0,
$$
a contradiction.  \qed

Lemma \ref{l1} implies $\cos y\ge 0$ on $[0,T]$.  Since $z(t)\in(0,{\pi\over 2}]$ and $\dot z(t)>0$ for $t\in(0,T]$, it follows that $\dot x$ is bounded below by a positive constant on $[\delta,T]$ for each $\delta>0$. Hence $X_a$ will reach $\partial D$ in finite time, and we denote by $t_a>0$ the first such positive time.
%$$
%t_a=\text{the first time $X_a$ reaches $\partial D$}.
%$$
The discussion before Lemma \ref{l1} shows that
$$
X_a(t_a)\notin  \{z=0\}\cup  \{x+y={\pi\over 2}\} \cup  \{y-x={3\pi\over 2}\},
$$
and hence 
$$
X_a(t_a)\in  \{x=0\} \cup \{z={\pi\over 2}\}.
$$

Let $S_0,S_1$ be the sets of all $a\in [0,{\pi\over 2})$ such that $X_a(t_a)\in \{x=0\}\setminus \{z={\pi\over 2}\}$ and  $X_a(t_a)\in \{z={\pi\over 2}\}\setminus \{x=0\}$, respectively.  We have $S_0\neq\emptyset\neq S_1$ since obviously $a\in S_1$ when $a$ is close to $\pi\over 2$, while Lemma \ref{l2} implies  $0\in S_0$.  Moreover, $\dot X_a(t_a)$ is transversal to $\partial D$ for any $a\in S_0\cup S_1$ (for $a\in S_0$ we have $\dot x(t_a)>0$ 
%is bounded below by a positive constant on $[{1\over 2} t_a,t_a]$ 
by the argument in the last paragraph, while for $a\in S_1$ obviously $\dot z(t_a)>0$).
It follows that $S_0,S_1$ are both relatively open in $[0,{\pi\over 2})$ and $t_a$ and $X_a(t_a)$ are continuous on them.
%it is easy to see that  $X_a(t_a)$ is not parallel to the boundary.  Hence $t_a$ is continuous with respect to $a$.   
Hence $[0,{\pi\over 2})\setminus(S_0\cup S_1)\neq\emptyset$, and then for any $\bar a$ from this set we must have
$$
X_{\bar a}(t_{\bar a})\in  \partial D\cap \{x=0\}\cap \{z={\pi\over 2} \} = \left\{(0,y,  {\pi\over 2}) \, \bigg| \,   y\in  [-{\pi\over 2},   {\pi\over 4}] \right\}.
$$

\medskip

\nit{\bf Step 2:}   We now use the symmetry of the ABC flow to show that for any $\bar a$ as above, $X_{\bar a}(t)=(x(t), y(t), z(t))$ is the desired solution.
%the orbit $X_{\bar a}(t)=(x(t), y(t), z(t)): \Rset\to \Rset^3$ is periodic along the $x$ direction.  In fact,  for 
For $t\in  \Rset$ let
$$
\tilde X(t)= (-\pi-x(-t),\ -y(-t),\ z(-t))  
$$
(this is the reflection across the line $(-\frac \pi 2,0)\times\Rset$) and
$$
\hat X(t)=\left (-x(2t_{\bar a}-t),\ y(2t_{\bar a}-t),\  \pi-z(2t_{\bar a}-t) \right)
$$
(this is the reflection across the line $\{x=0\}\cap\{ z=\frac\pi 2\}$).
Clearly,  both $\tilde X$ and $\hat X$ are solutions to the 1-1-1 ABC flow system.  Since  $\tilde X(0)=X_{\bar a}(0)$ and $\hat X(t_{\bar a})=X_{{\bar a}}(t_{\bar a})$, we have $X_{\bar a}=\tilde X=\hat X$.   Thus
$$
(x(-t_{\bar a}), \  y(-t_{\bar a}),\ z(-t_{\bar a}))=X_{\bar a}(-t_{\bar a}) =\tilde X(-t_{\bar a})= \left(-\pi,-y(t_{\bar a}),\frac \pi 2 \right), 
$$
and then  
$$
\hat X(3t_{\bar a})=(-x(-t_{\bar a}), y(-t_{\bar a}),  \pi-z(-t_{\bar a}))= \left (\pi ,-y(t_{\bar a}), \frac\pi 2 \right).
$$
So
$$
X_{\bar a}(3t_{\bar a}) =\hat X(3t_{\bar a})= X_{\bar a}(-t_{\bar a})+(2\pi,0,0),
$$
and it follows from $2\pi$ periodicity of the ABC flow that 
$$
X_{\bar a}(t+4t_{\bar a})=X_{\bar a}(t)+(2\pi,0,0).
$$
for each $t\in\Rset$.
\qed

\bibliographystyle{plain}

\end{document}